\documentclass[review]{elsarticle}

\usepackage{graphics}
\usepackage{subfigure}

\usepackage{amssymb}

\makeatletter
\def\ps@pprintTitle{%
 \let\@oddhead\@empty
 \let\@evenhead\@empty
 \def\@oddfoot{}%
 \let\@evenfoot\@oddfoot}
\makeatother

\begin{document}

\begin{frontmatter}



\title{On numerical approaches to the analysis of topology 
of the phase space for dynamical integrability}


\author{Vladimir Salnikov}

\address{ Laboratoire de Math\'ematiques de l'INSA de Rouen, \\ 
Avenue de l'Universit\'e 76801 Saint-\'Etienne-du-Rouvray Cedex, France \\
              Tel.: +33 2 32 95 65 52\\
              vladimir.salnikov@insa-rouen.fr }
              
\begin{abstract}
In this paper we consider the possibility to use numerical simulations 
for a computer assisted qualitative analysis of dynamical systems. 
We formulate a rather general method of recovering the obstructions 
to dynamical integrability for the systems that after reduction have a
small number of degrees of freedom. 
We generalize this method using the results of KAM theory and stochastic 
approaches to the families of parameter depending systems. 
This permits to localize possible integrability regions in the parameter space. 
We give some examples of application of this approach to dynamical systems
having a mechanical origin.
\end{abstract}

\begin{keyword}
Dynamical integrability \sep numerical approach \sep
phase space topology \sep KAM theory \sep Monte-Carlo method \sep pendulum-type systems
\end{keyword}

\end{frontmatter}


\section{\label{sec:intro}Introduction}
This paper is a part of a series of works devoted to description of the 
possibilities of application of numerical methods to qualitative analysis 
of dynamical systems. 
The main subject that we are going to discuss will be dynamical integrability.
The problem of integrability has been studied since the middle of the XIX century,
when the question of primary interest was to be able to integrate the 
system of differential equations by inversion of functions and quadratures, that is to give 
a more or less explicit solution. Nowadays with the development of 
methods relating symplectic geometry and dynamics the notion is extended, 
namely one studies the existence of an appropriate number of conserved quantities
(first integrals, invariant measure, etc), possessing some properties.
One of the main motivations to study integrable systems is related to the fact that they 
are in a sense more `regular' than the generic ones: for such systems the questions of perturbation and stability are of particular interest. But this does not mean that 
non-integrable systems are `bad', since for them one can expect interesting non-linear 
behaviour, chaotization etc. Thus, there are two natural directions of the study of 
integrability: one is a search for non-artificial examples of integrable 
systems, the other is a rigorous proof of non-integrability for given dynamical 
systems -- we will address both of them.

In this paper we will consider the  integrability in the Liouville--Arnold sense, namely 
for a hamiltonian system with $n$ degrees of freedom it is the existence of $n$ independent 
first integrals in involution (\cite{AKN}). In particular we will pay attention to 
the independence condition from the topological point of view.
We suggest a constructive method of a computer assisted analysis of integrability 
of systems with small dimensional configuration space which is based on the 
study of the topology of the phase space. 
We generalize it using the results of the Kolmogorov-Arnold-Moser theory to systems with parameters. 
The logic of the method is explained in details via an example of 
pendulum-type systems, then comments on the range of its applicability are presented. 
We also mention other examples of application of this method to some 
concrete dynamical systems having a mechanical origin.

\section{Method of sections} \label{sec:sec}
Consider an autonomous hamiltonian system with a two-dimensional configuration 
space $Q$, its phase space $T^*Q$ is of dimension $4$. 
Since the right hand sides of the equations of motion do not depend 
explicitly on time, any trajectory of this systems belongs to a constant hamiltonian (energy)
level hypersurface having the dimension $3$. For complete integrability 
of the system with two degrees of freedom
the existence of another first integral independent with energy integral is needed -- 
in this case any trajectory would belong to the $2$-dimension manifold 
defined by the intersection of the level surfaces of these first integrals.

\subsection{Idea of the method}
This simple topological consideration alone does not 
give a method of analysis of integrability, but together with a good visualization 
algorithm it permits to give an answer for a rather large class of systems.
Effective visualization of the dynamics in a $4$-dimensional phase space is not an easy 
task, therefore we consider the intersection of a trajectory of the 
system with $2$-dimensional planes in it. Let us discuss the possible results 
of this intersection. 

The dimension of intersection of two generic manifolds 
of dimensions $n_1$ and $n_2$ in an 
$N$-dimensional space is given (\cite{DNF}) by the equation 
\begin{equation} \label{dims}
  dim = n_1 + n_2 - N.
\end{equation}
In our case the manifold swept by the trajectory intersects a $2$-dimensional plane 
($n_2=2$) in the phase space of dimension $N=4$.
Energy conservation guarantees that the dimension ($n_1$) of this manifold is at most $3$.
Thus we can observe two possible cases depending on the existence of additional first integral:
either $n_1 = 2$ and the intersection is $0$-dimensional (finite set of points), 
or $n_1=3$ and the intersection is of dimension $1$ 
(finite set of curves). 
It is clear that the presence of curves in the intersection is an obstruction to 
complete integrability. Absence of curves however does not directly mean 
integrability, since it can be the consequence of various reasons -- we will discuss 
them more precisely in  section \ref{sec:kam} in the context of generalization 
of this method.

Let us be more explicit on this idea about intersection and prove 
in a particular case the formula (\ref{dims}). 
Choose  canonical coordinates $q_1, p_1, q_2, p_2$ of the phase space $T^*Q$
(assume for simplicity that it can be done globally).
Then an arbitrary two-dimensional plane (or, better to say two-dimensional linear subspace)
is given by the system of equations
\begin{eqnarray} \label{planes}
    a_1 q_1 + b_1 p_1 + c_1 q_2 + d_1 p_2 = e_1, \nonumber \\
    a_2 q_1 + b_2 p_1 + c_2 q_2 + d_2 p_2 = e_2,
\end{eqnarray}
for some constants $a_i, b_i, c_i, d_i, e_i, i = 1,2.$
Since the system is autonomous any trajectory belongs to the energy level
\begin{equation} \label{energy}
    I_1 \equiv H = h,
\end{equation}
where $H$ is the hamiltonian function defining the dynamics of the system.
The conditions (\ref{planes}, \ref{energy}) on $4$ coordinates in the phase space 
define a $1$-dimensional manifold. If there exists another first integral
\begin{equation} \label{2int}
  I_2 = const
\end{equation}
the system (\ref{planes}, \ref{energy}, \ref{2int}), if it is compatible, 
admits a finite set of solutions.  

It is now easy to see that it is important to pay attention to the independence 
of the conditions (\ref{planes}) from (\ref{energy}) and  (\ref{2int}).
This is not difficult to guarantee for the energy integral, since it is non-linear 
and usually known; if the other first integral exists and the equation (\ref{2int}) 
turns out to be dependent with (\ref{planes}) it means that the system admits some
degeneracy, that will be seen in the numerical simulation.
To avoid ``false detection'' by this method it is enough to 
consider the intersection with several mutually independent planes, that is 
done in practice\footnote{I am thankful to the anonymous referee for underlining
the importance of this remark.}. 
Namely, one can modify freely the coefficients $a_i, b_i, c_i, d_i$,
that corresponds to ``rotations'' of the planes, as well as $e_i$ that corresponds to 
considering the parallel planes, -- for a generic situation both 
operations should not modify drastically the topology of the intersection. 
In case of doubt it may be also useful to check if for a fixed choice of 
coordinates on the plane (\ref{planes}) there exists an integral polynomial 
in the complementary coordinates in the whole phase space, or at least 
if the intersecting plane is invariant under the action of the flow of the system.

\subsection{Triple pendulum}
As one of the main examples of application of the method in this work we will consider 
the free flat motion of pendulum-type systems.
A multiple pendulum is the system of point masses, connected by weightless 
inextensible rods, the first of the points being fixed.
For a triple pendulum these conditions read
$$
({\bf r}_i - {\bf r}_{i-1})^2 = l_i \qquad i = 1, 2, 3,
$$
where ${\bf r}_i$ corresponds to the $i$-th mass, the fixed point is $ {\bf r}_0$. 
One can also consider a more general case of constraints given by 
arbitrary  polynomials of degree $2$, which however does not always correspond to 
a physical configuration.

The choice of this system is motivated by several factors. 
Let us note that the case of a double pendulum is well studied. 
In particular its free motion on the plane is a classical example of a completely 
integrable system (see for example \cite{whittaker}).
In the above notations taking the 
angular momentum for $I_2$ one obtains two independent conditions (\ref{energy}, \ref{2int}),
having a $0$-dimensional intersection with (\ref{planes}).
But already in the presence of gravity  some chaotic behaviour has been observed, 
that is the trajectory 
is rather dense on the energy level-surface, that can be seen on the sections.
A rigorous proof of non-integrability of this problem is however 
a subtle question (\cite{ramis_question}).
The problem of control has also been studied for this system (\cite{chern}),
among interesting results one can mention that a double pendulum can be stabilized in the 
upright position by controlling only one degree of freedom.

Meanwhile the dynamics of the triple pendulum is almost not studied, although 
it is interesting for applications. 
We have shown rather directly (without going much into details about topology, 
but writing explicitly the coordinates) that the behaviour of the system 
is rather irregular and it is not integrable\footnote{This work has been carried out in collaboration 
with V.L.~Golo and resulted in the first instance of application of the methods of sections}.
We will not describe this example in full details since we have already done 
it partially in  \cite{imec} and in \cite{cirm} and the result will also follow as a particular 
case from the section \ref{sec:kam}. Let us just note that 
in the original setting the system obviously has $3$ degrees of freedom, 
but since it admits the angular momentum first integral 
it can be reduced by the Routh transform (\cite{theormech}) to the one having the $4$-dimensional 
phase space, that is precisely in the range of application of the method of sections. 
Let us also note, that this reduction is profitable to understand the topology 
of the phase space, but to perform numerical simulations
it is more convenient to use the Lagrange multipliers (\cite{lagrange})
and then compute the reduced coordinates from the cartesian ones.

\begin{center}
 \begin{figure}[htp] \centering
     \includegraphics*[height=4in]{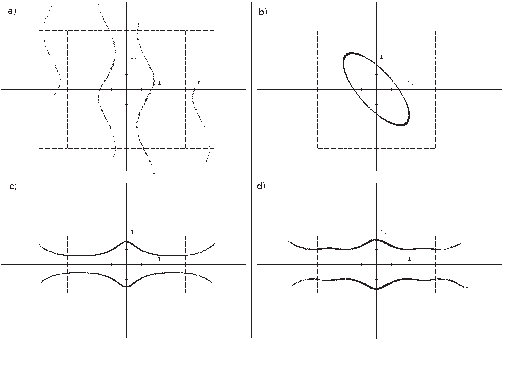}
      \caption{
            \label{3pen_sections}
            Intersection of the trajectory of the triple pendulum with the planes: \newline
             a). $(\dot\beta_1 = 1, \dot \beta_2 = 1)$, 
	     b). $(\beta_1, \beta_2)$,
	     c). $(\beta_1, \dot \beta_1)$,
	     d). $(\beta_2, \dot \beta_2)$. 
            }
 \end{figure}
\end{center}

A typical result of the numerical simulations is shown on figure \ref{3pen_sections}.
The coordinates in the reduced phase space are the angles $\beta_1, \beta_2$
between the segments of the pendulum and their derivatives. 
In the hamiltonian formalism one should actually use momenta instead of 
velocities, but we use the natural duality of $T^*Q$ and $TQ$, so it gives an equivalent picture 
from the topological 
point of view. 
The intersecting planes are chosen to be parallel to the coordinate ones. The presence of
curves in the intersection clearly shows the non-integrability.

\subsection{Satellite dynamics}
\begin{center}
 \begin{figure}[htp] \centering
    \includegraphics*[width=0.9\linewidth]{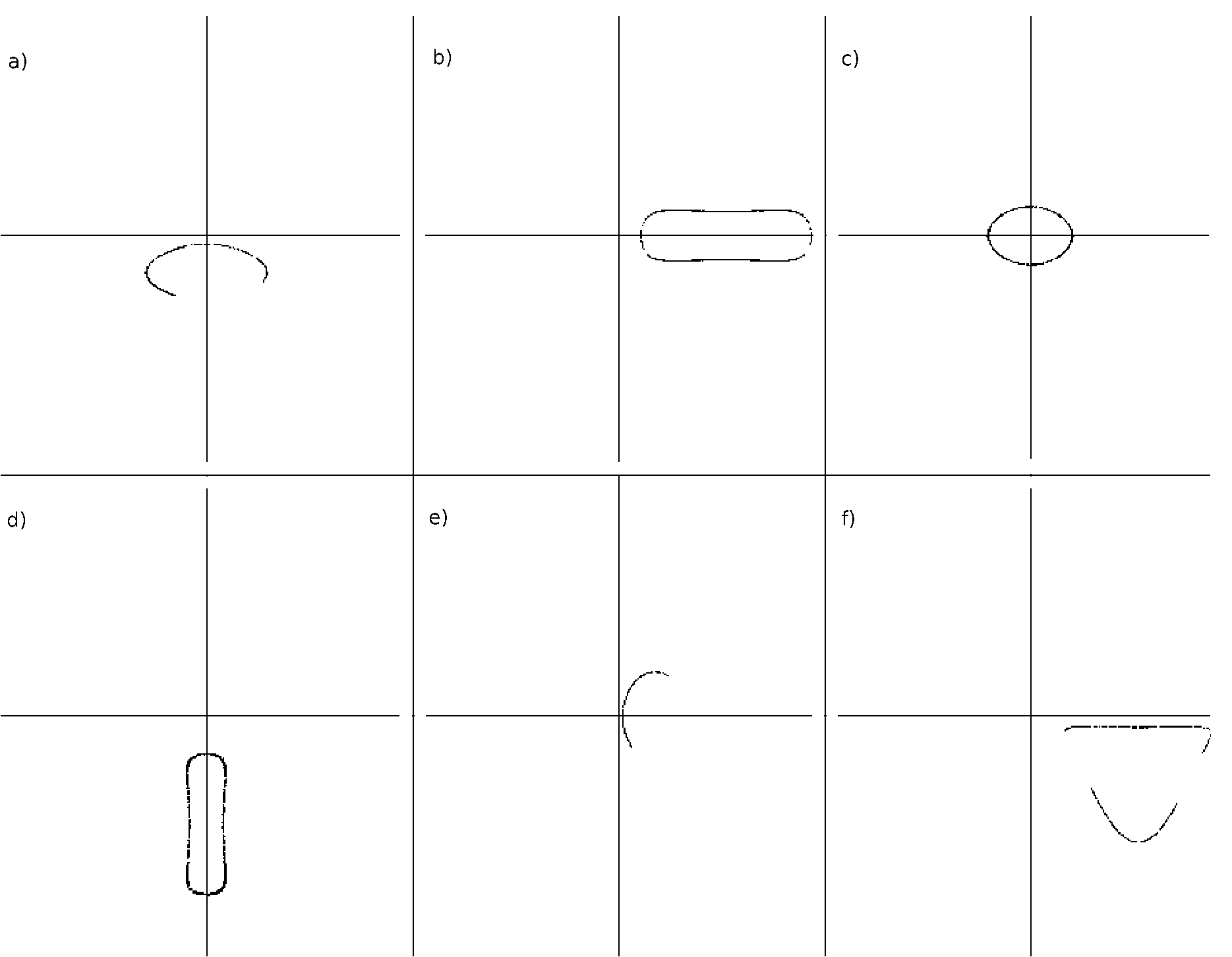}
    \caption{
              \label{sat}
               Symmetric satellite: intersection of the trajectory 
               with the planes, parallel to the coordinate ones, containing the 
               point $\psi = 0.28, \theta = 0.82, p_{\psi} = 0.15, p_{\theta} = 0.37$.   
         }
 \end{figure}
\end{center}

Let us consider another mechanical example for which the method is applicable:
the motion of the dynamically symmetric satellite on a circular orbit.
In the orbital coordinate system, the axes of which are directed along the 
radius of the orbit of the center of masses of the satellite, its normal and its binormal, 
the position of the satellite is given by the three Euler angles
$\psi, \theta, \varphi$. Denoting the corresponding momenta by $p_\psi, p_\theta, p_\varphi$,
one obtains the following hamiltonian function:
\begin{eqnarray}
  H = \frac{p^2_\psi}{2 \sin^2\theta} + \frac{p^2_\theta}{2} - p_\psi ctg(\theta)\cos(\psi) - 
  \alpha \beta p_\psi \frac{\cos \theta}{\sin^2\theta} - p_\theta \sin(\psi) +  \nonumber \\
   + \alpha \beta \frac{\cos \psi}{\sin \theta} + \frac{\alpha^2 \beta^2}{2\sin^2 \theta} + \frac{3}{2}(\alpha - 1)\cos^2\theta, 
\label{sat_ham}
\end{eqnarray}
where $\alpha = C/A$; $A,B,C$ are the principal moments of inertia ($A = B$).
The coordinate  $\varphi$ is cyclic, thus we can fix the corresponding momenta 
$p_\varphi = \alpha\beta = const$, where $\beta$ denotes the ratio of the 
orbital angular velocity and the projection of the angular velocity of the satellite to 
its symmetry axis. Let us consider the values of  $\alpha = 4/3, \beta = 0$, 
for which the stability of equilibrium solutions is studied in 
\cite{bardin}. 
The reduced system has two degrees of freedom and is described 
(up to redefinition of coordinates) by the hamiltonian:
$$
  H = \frac{\displaystyle p_{\psi}^2}{\displaystyle 2 \sin^2\theta} + \frac{\displaystyle p_{\theta}^2}{\displaystyle 2}
   - p_{\psi} + \frac{1}{2} \sin^2\psi \sin^2 \theta.
$$
That is we can again apply the method of sections.
The intersections that one obtains (figure \ref{sat}) show its non-integrability.
Let us however note, that to obtain the sections containing curves for this problem 
we had to analyze a lot of different initial conditions. 
This indicates that the system can be \textit{locally integrable}, 
that is possess an additional first integral for a certain subset of initial data.

\subsection{Details of the method of sections}
Let us now discuss the possibilities of application of the method of sections 
to other systems, in particular let us comment on the 
class of systems for which it can result in the precise conclusion 
on integrability.
As we have noted in the very beginning the natural restriction comes from the 
dimension of the phase space, or more precisely from the 
possibility to reduce it. This however covers a large class of systems for which  
integrability is an open question, such as  dynamics of a triple lattice, 
the motion of a mass point in a symmetric potential in a $3$-dimensional space, 
dynamics of geodesics on curved surfaces etc.
Let us note that the examples above represent a typical situation 
in the integrability analysis, when for a system with $n$
degrees of freedom $(n-1)$ first integrals are known and one is interested in the
existence of a supplementary one. Then an important feature of the method is that 
the explicit reduction of the system is not needed, that is one can 
study the trajectory of the initial system and only use the fact that the 
reduction can be done. 

The only difficulty is the choice of convenient coordinates 
for constructing the sections. But this problem can be naturally solved when the system 
admits the angle-parametrization like in the described examples. Note that such systems 
often arise in the applications. According to the theorem by V.V.~Kozlov (\cite{kozlov_sym})
it is interesting to consider the systems with the genus of the configuration space 
at most $1$, since if it is not the case one knows that the system is not 
analytically integrable; but it means that basically any choice of coordinates 
has the structure of angles. We can however study more general surfaces 
if we are interested in a larger class of first integrals. Then it is good to 
make sure that the motion takes place in a bounded domain of the phase space, 
to accumulate a sufficiently representative intersection -- 
this condition is also rather natural at least for the systems with compact 
configuration space. 

Let us also note a technical difficulty of the method. 
To conclude non-integrability one needs to compute a rather long trajectory, 
since the intersection of it with a surface of small dimension is a rear event
in contrast, say,  to Poincar\'e sections where every loop gives a point in the plane.
But this reduction of the dimension is done intensionally, as it permits to 
classify different cases more effectively: it is much easier to 
distinguish points from curves than curves from thin domains that can occur in the 
Poincar\'e sections -- this process can even be done automatically without 
interaction with the user.
Our method also permits not to pay much attention to 
transversality of the trajectory to the planes. All this only results in the need to 
apply reliable algorithms of numerical integration.

Another thing that we haven't yet addressed in details is the interpretation 
of the ``empty'' ($0$-dimensional) intersections. 
The existence of such sections does not permit to make any direct conclusion since it only means that 
a given trajectory is by chance more regular.
But in the next section we will explain the origin of such sections for non-integrable 
systems and discuss the possibility to use them to qualitatively describe the 
behaviour of the systems.

\section{Generalization of the method of sections via the KAM theory} \label{sec:kam}
The method of sections as we have seen permits to prove non-integrability 
for a given dynamical system. This problem however is rather special, more often
one wants to find the relation between the parameters when integrability is 
possible. In this section we will consider an example of a two parameter 
system showing how one can extend the method of section to 
this problem. Despite the discrete nature of the method it permits 
(via some extra mathematical considerations) also to draw conclusions 
for a continuous range of parameters.

\subsection{Pendulum-type systems} 
\begin{center}
 \begin{figure}[htp] \centering
    \includegraphics*[height=1.3in]{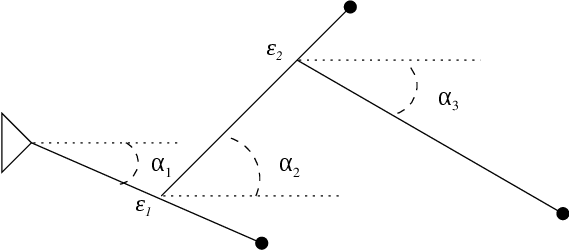}
    \caption{              \label{gen3pen}
               A pendulum type system parametrized by the angles. 
            }
 \end{figure}
\end{center}

Consider the free flat motion of a system obtained from the triple pendulum 
by moving the fixing points of the next segment along the previous one 
(figure \ref{gen3pen}). We shall call such objects pendulum-type systems.
In particular they describe the motion of a physical pendulum, 
i.e. rigid bodies fixed between themselves.

The configuration of this system is described by two-dimensional vectors 
$\mathbf{r_i}$ and velocities $\mathbf{v_i}$ ($i=1, 2, 3$) with the following constraints:
 \begin{eqnarray}
\varphi_1 = (\mathbf{r_1} - \mathbf{r_0})^2 - l_1^2 = 0, \nonumber \\
\varphi_2 = (\mathbf{r_2} - \varepsilon_1\mathbf{r_1})^2 - l_2^2 = 0, \nonumber \\
\varphi_3 = (\mathbf{r_3} -  ( \varepsilon_1\mathbf{r_1} + \varepsilon_2( \mathbf{r_2} -  \varepsilon_1 \mathbf{r_1} ) ) )^2 - l_3^2 = 0,
  \end{eqnarray}
the point $ \mathbf{r_0}$ is fixed. 
The parameters $\varepsilon_i$ characterize linearly the fixing point 
of the $(i+1)$-st segment to the $i$-th one: $\varepsilon_i = 1$ corresponds to 
fixing at the endpoint (like in the multiple pendulum), 
$\varepsilon_i = 0$ -- to the fixed point (like two non-interacting pendula with the 
same fixed point).

Note that for $\varepsilon_1 = \varepsilon_2 = 1$ we obtain the non-integrable case 
of the triple pendulum studied before. 
For $\varepsilon_1 = 0$ and  any $\varepsilon_2$ 
the system decouples to non-interacting simple pendulum and a 
double pendulum-type system, and becomes obviously integrable.
And in the particular case when both $\varepsilon_1 = \varepsilon_2 = 0$, 
it decouples to three non-interacting simple pendula.
For other values of  $\varepsilon_i$ we will be again interested in the existence of 
a sufficient number of independent first integrals.

\subsection{Topology of the phase space}
For all values of $\varepsilon_1, \varepsilon_2$ the system is invariant 
under rotation around the fixed point  $ \mathbf{r_0}$, 
thus it possesses the Noether integral of angular momentum.
Since the position is still characterized by the three angles 
(figure \ref{gen3pen}), this integral permits to perform Routh reduction
to the $4$-dimensional phase space. To make this reduction explicitly
we again study the dynamics in terms of the angles between the segments of the 
pendulum-type system $\beta_i = \alpha_{i+1} - \alpha_i, \quad i = 1, 2$ 
and their derivatives.
Now for any fixed couple of  $(\varepsilon_1, \varepsilon_2$), 
 we can apply the method of sections.

\begin{figure}[ht]
\centering
\subfigure[\, $\varepsilon_1 = 0.7, \varepsilon_2 = 0.8$]{
    \includegraphics*[width=2.in]{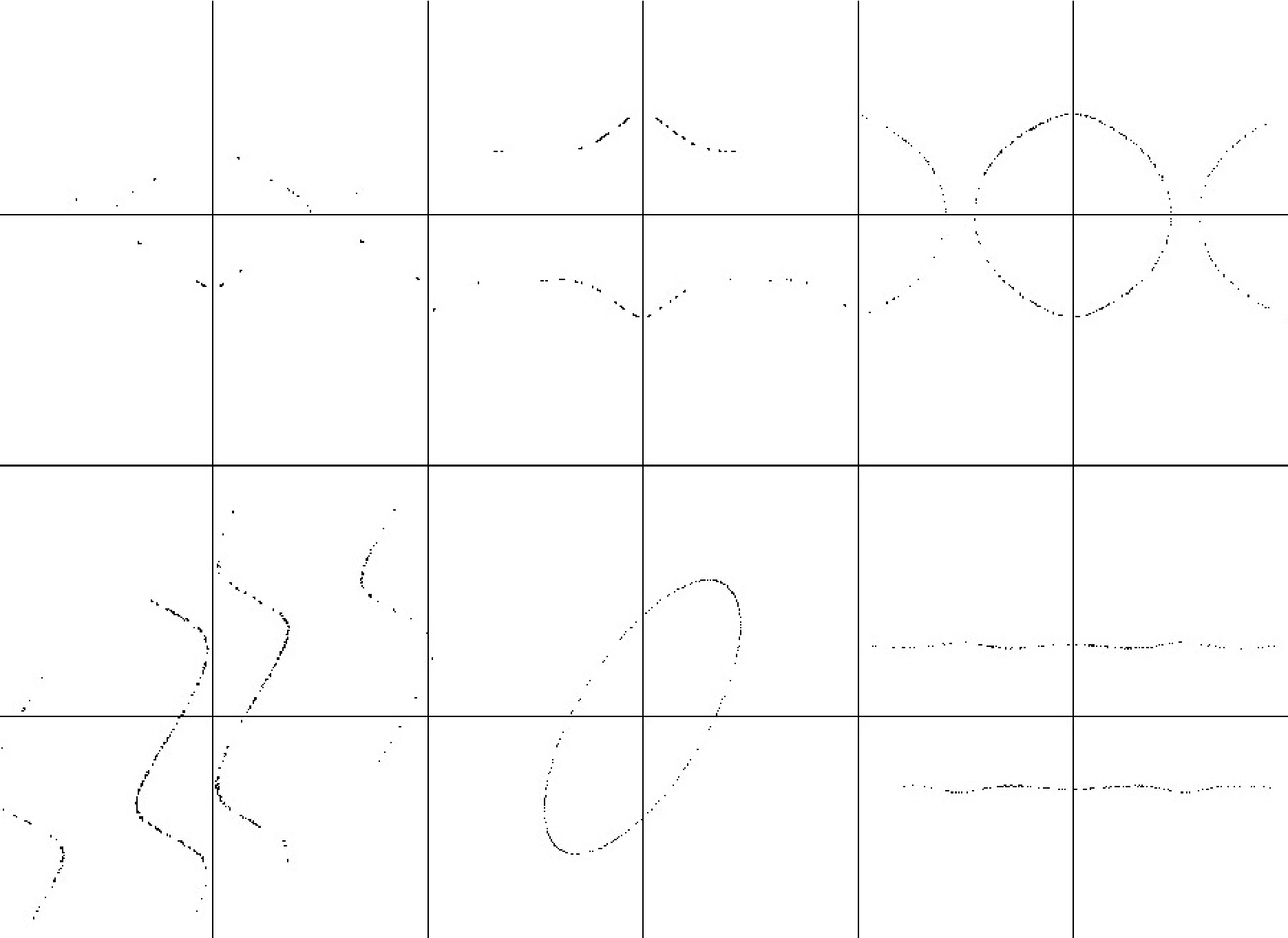} \label{7080}
 }
\subfigure[\, $\varepsilon_1 = 0.72, \varepsilon_2 = 0.74$]{
    \includegraphics*[width=2.in]{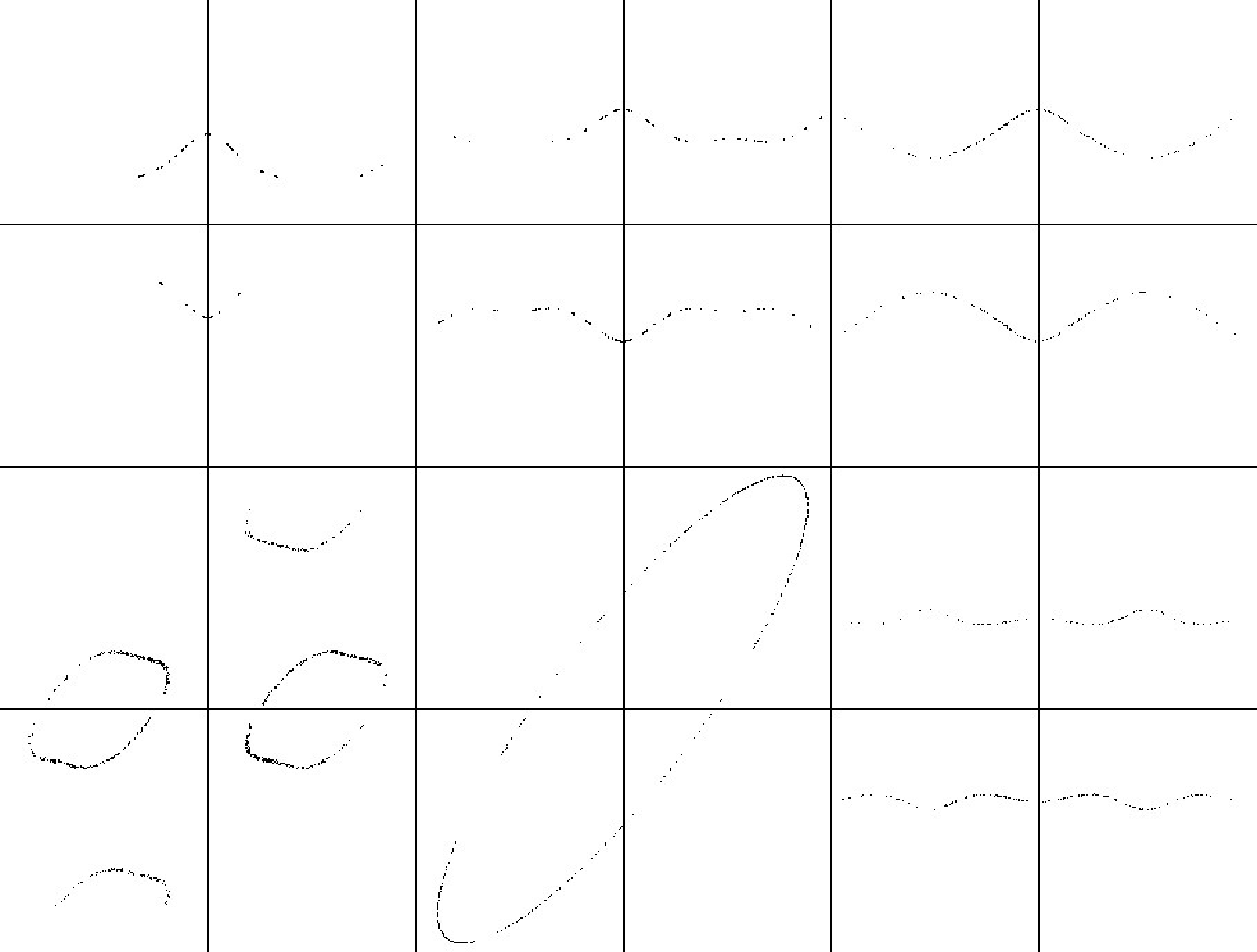} \label{7274}
 }
 
\subfigure[\, $\varepsilon_1 = 0.72, \varepsilon_2 = 0.76$]{
    \includegraphics*[width=2.in]{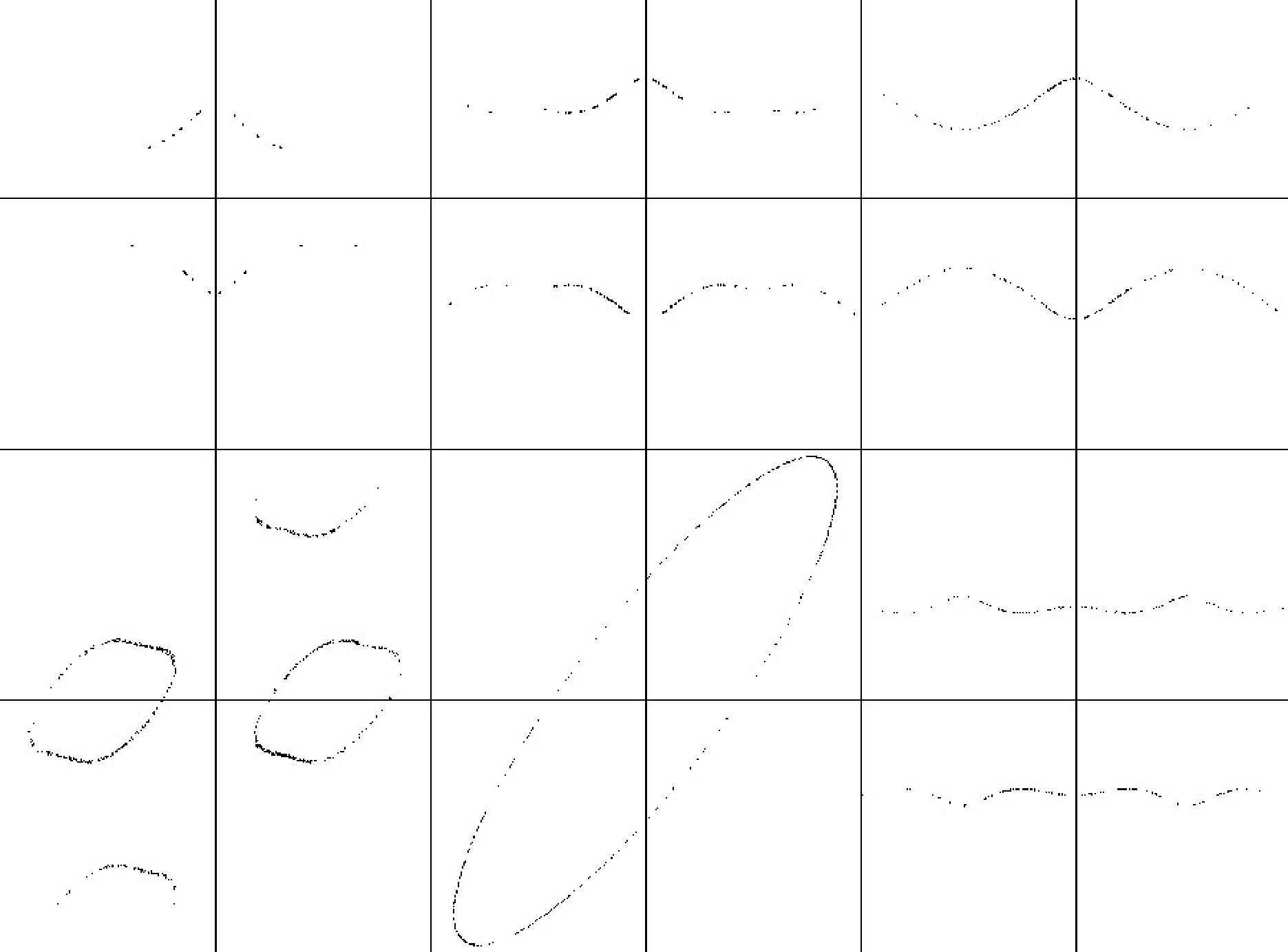} \label{7276}
 }
\subfigure[\, $\varepsilon_1 = 0.74, \varepsilon_2 = 0.78$]{
    \includegraphics*[width=2.in]{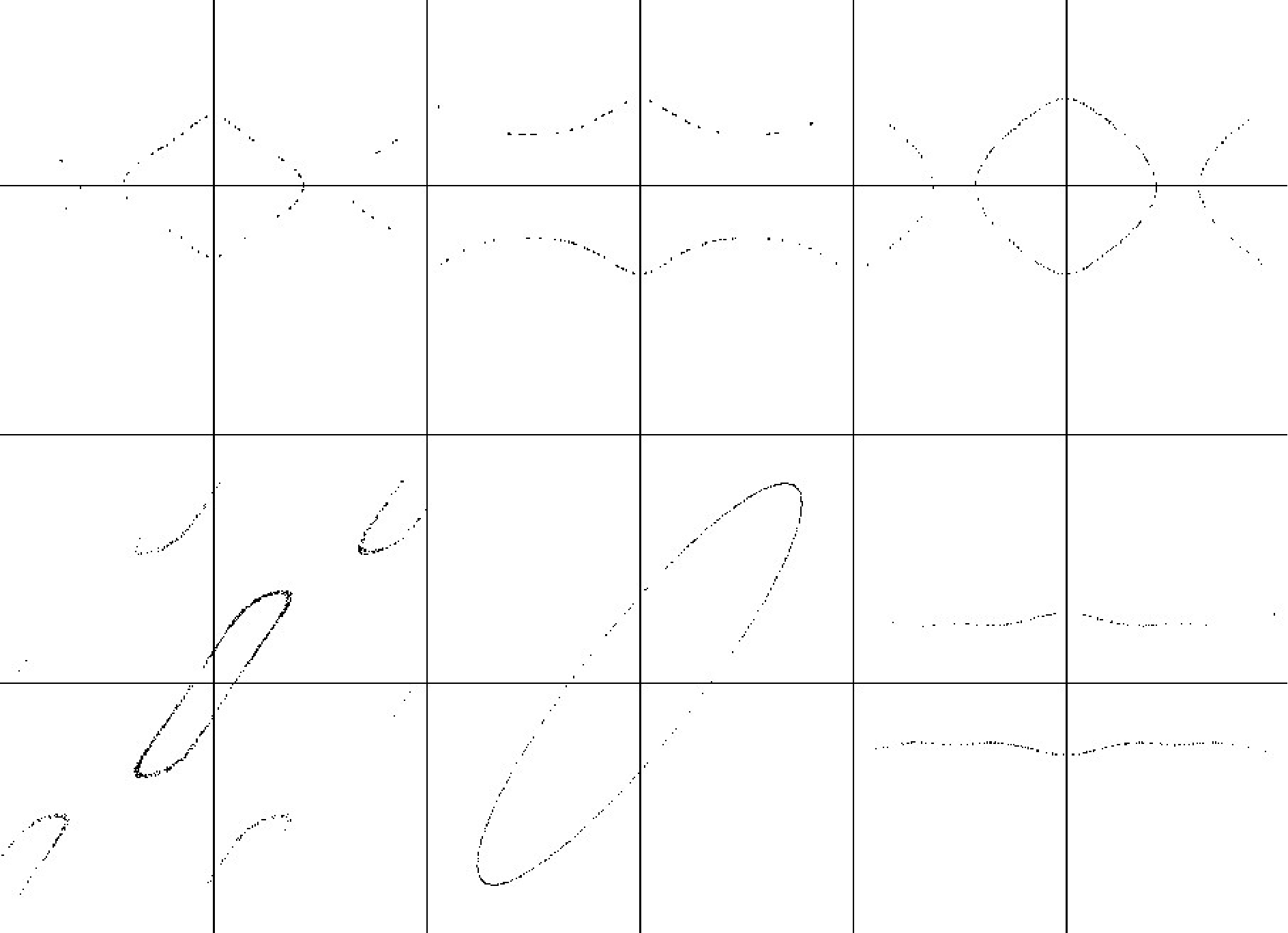} \label{7478}
 }
\label{4figs} 
\caption{Sections for various values of $\varepsilon_1, \varepsilon_2$}
\end{figure}

Typical results for the pendulum-type systems are represented on figures \ref{7080}--\ref{7478}, each
showing the intersection with the planes 
$(\dot \beta_2 = 0, \beta_2 = 0), (\dot \beta_2 = 0, \beta_2 = 0), (\dot \beta_1 = 0, \beta_2 = 0), 
(\dot \beta_2 = 1, \dot \beta_1 = 1), (\beta_2 = 0, \beta_1 = 0), (\dot \beta_1 = 0, \beta_1 = 0)$. 

This figures show that the sections contain curves, it means the corresponding systems 
do not admit an additional first integral. The similar figures have been obtained 
(i.e. initial conditions found\footnote{A sufficiently time-consuming numerical experiment, performed 
on the cluster of ICJ, Universit\'e Claude Bernard Lyon 1.}) for all couples 
$\varepsilon_1, \varepsilon_2 \in [0,1]\times[0,1]$ with the step $0.02$, 
except the trivial ones ($\varepsilon_1 = 0$), that means non-integrability of 
the corresponding systems.

\subsection{Generalization to continuous parameters}
Note that the sections for close values of parameters and close initial data
look rather similar (figures \ref{7274}, \ref{7276}). 
One thus is tempted to conjecture that 
the systems between  two values of parameters are also non-integrable.
Alone this statement is certainly false, but we can make sense out of it 
by some extra consideration. Let us make it more precise.

The results of Poincar\'e (\cite{p1,p2}) 
show that the systems obtained from an integrable one by a perturbation 
usually fail to be integrable. However according to the 
Kolmogorov-Arnold-Moser theorem (see for example \cite{AKN}) 
for the systems close to integrable a set of positive measure of invariant tori
preserves its topology. In our terms it means, that when the initial point 
even of a non-integrable system lies on such a torus the sections 
do not contain curves. And when we approach an integrable system in the space of
parameters the probability to be on such a torus increases.

\begin{figure}[ht]
\centering
\subfigure[\, Depending on  $(\varepsilon_1,  \varepsilon_2) \in  (0,1)\times(0,1) $. \newline
	       Scale: from $0$ -- white, to $1$ -- black. 
	       ]
	       {	\includegraphics*[width=0.6\linewidth]{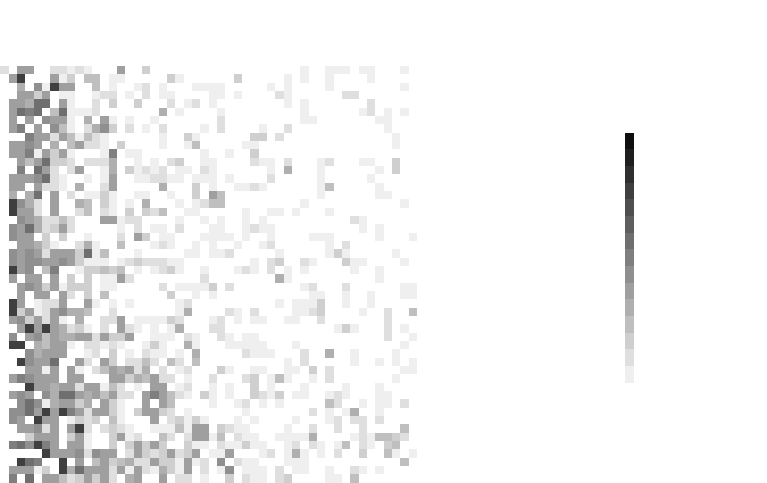}       \label{stat}}
\subfigure [\, Depending only on \newline $\varepsilon_1 \in (0,1)$]  
{   \includegraphics*[width=0.3\linewidth]{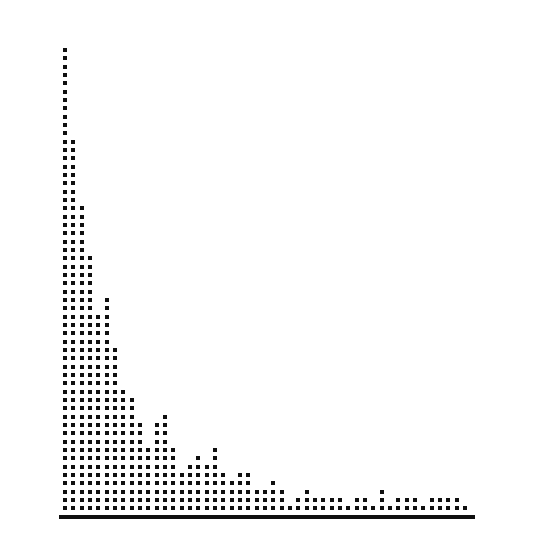}   \label{stat2}}
\label{2figs} 
\caption{Distributions of the Liouville tori preserving their topology. \newline
The maxima are accumulated along $\varepsilon_1 = 0$}
\end{figure}

This effect is indeed observed in the neighborhood of the integrable case
($\varepsilon_1 = 0$) for arbitrary chosen initial points.
The figure \ref{stat} shows the distribution obtained by the Monte-Carlo method (\cite{monte-carlo})
of such tori, 
i.e. the proportion 
of ``empty'' sections if one starts from a pseudo-random point in the phase space 
for the given values of parameters $(\varepsilon_1, \varepsilon_2)$.
The figure \ref{stat2} shows the same distribution only depending 
on $\varepsilon_1$. One sees from them, that except the segment 
$\varepsilon_1 = 0$, in the square $(\varepsilon_1, \varepsilon_2) \in [0,1]\times[0,1]$ 
there are no pronounced maxima, that shows non-integrability of all the systems 
with parameters different from $\varepsilon_1 = 0$.
That is in the family of triple pendulum-type systems the only integrable situation 
corresponds to decoupling the systems to non-interacting subsystems of smaller dimensions.
Let us note that this conclusion is coherent with the above discussion 
about the triple pendulum as well as with the recent result of \cite{flail}
concerning the ``flail'' triple pendulum.

\textbf{Generalization of the method of sections.}
It is clear that the same approach is not limited to the above example. 
We can apply it to any family of systems with arbitrary number of parameters 
the phase space of which can be uniformly reduced to the $4$-dimensional one.
One needs to choose a domain $D$ in the parameter space, 
apply the method of sections to a sufficient number of 
pseudo-random initial points for all values of parameters from a rather 
dense set in $D$ and compute the estimation of the probability 
to start from an invariant torus preserving its topology.
The subdomains of $D$ corresponding to maxima of the obtained distribution 
are the candidates for integrability.

Let us stress again that the idea to use the method of sections together with the results of the
KAM theory is crucial for applications: instead of studying the property of a concrete system to be integrable we analyze the qualitative behaviour of a family of systems. 
We are thus replacing a point-wise property in the space of parameters 
which can be easily missed in numerical simulations by a local property that
is much better observable. And the only effective way to study such local 
properties is to have an easily (automatically) distinguishable characteristic 
like the one we suggest in the method. 
It is important also to note that in contrast to the 
visualization algorithm from the section \ref{sec:sec}, the key idea of the method --
accumulation of ``empty'' sections around distinguished points in 
the parameter space -- can be extended to higher dimension of the 
(reduced) phase space, up to a classification algorithm based on the equation 
(\ref{dims}) that 
becomes technically more involved.

Continuing the remark from the previous section let us note that one should not 
disregard the time needed to perform the computation related to application 
of the generalized method of sections. We are studying a sufficiently large 
number of trajectories of the system for each set of values of the parameters, 
but these simulations are independent from each other. 
And since for any of them the decision can made automatically, 
the process is perfectly parallelizable and very well suitable for the usage 
of the modern multiprocessor computing systems based for example on the GPU.

\section{Conclusion}
In this paper we have described the topological properties of the phase space of
the small dimensional integrable dynamical systems in the form convenient for 
verification. We have proposed a method of computer assisted 
analysis of integrability for a given dynamical system and generalized 
it to the continuous range of parameters.
We have applied it to several systems the dynamics of which is interesting in the 
context of stability problems in mechanics.

This approach is a natural application of the classical Monte-Carlo method,
permitting by a well-developed technique, using the pseudo-random quantities to
construct a qualitative picture of the behaviour of a deterministic system.
Together with the results of the Kolmogorov-Arnold-Moser theory it 
permits to answer rather complicated questions on integrability.
Let us also note that on top of the described applications we were able 
to recover some know results from the theory of dynamical systems, like the 
classical integrable cases for the problem of the motion of a heavy rigid body 
with a fixed point (\cite{matapli}). 
Some work in progress also concerns qualitative analysis of 
dynamical systems not directly related to integrability: 
since the presented method permits in a sense to characterize 
quantitatively the chaotization of the system it helps to 
``measure'' how generic are regular/irregular trajectories.

As we have mentioned in the beginning of this paper there are other ways 
to study the qualitative behaviour of dynamical systems using numerical methods. 
For example the meromorphic non-integrability of the examples discussed in the section
\ref{sec:sec} can be proved by constructing the monodromy group (\cite{Ziglin_num})
and application of the results of \cite{Ziglin}.
Let us also briefly mention that one of the main motivations to develop these methods 
for us is studying more involved problems like systems with delay naturally appearing in 
biological modeling or relativistic celestial mechanics. There an important question 
is to study the existence of bounded solutions or the solutions of the limit cycle 
type. In the above language they will correspond to more regular trajectories than the 
arbitrary ones, and can be therefore localized by the likewise methods.\newline

\textbf{Acknowledgements.}
This work has been mostly done in the Claude Bernard Lyon 1 University and 
the Dorodnitsyn Computing Center of Russian 
Academy of Sciences -- the author is thankful to Sergey Stepanov for 
constant attention. The author also thanks Boris Bardin for the suggestion to look at 
the satellite dynamics 
within the framework of the developed approach as well as for providing 
the manuscript of the article \cite{bardin} in preparation, 
and useful discussions of the obtained results.\\


%






\end{document}